\input ./style/arxiv-general.cfg
\documentclass[MSNbibl,number,citesort,seceqn,dvips]{arxbj}
\makeatletter
   \@ifpackageloaded{graphicx}{}{\usepackage{graphicx}}
\makeatother
\usepackage{upgreek}

\volume{21}
\issue{4}
\pubyear{2015}
\firstpage{2569}
\lastpage{2594}
\doi{10.3150/14-BEJ655} 
\docsubty{FLA}

\makeatletter

\newcommand{\rrvert}{\vert}
\newcommand{\llvert}{\vert}
\newcommand{\eqref}[1]{(\ref{#1})}
\newtheorem{theorem}{Theorem}[section]
\newtheorem{proposition}[theorem]{Proposition}
\newtheorem{lemma}[theorem]{Lemma}
\newproclaim{definition}[theorem]{Definition}
\newremark{remark}{Remark}[section]

\makeatother

\begin{document}
\begin{frontmatter}

\title{Efficient pointwise estimation based on discrete data
in ergodic nonparametric diffusions}
\runtitle{Efficient discrete data estimation in diffusions}

\begin{aug}
\author[A]{\inits{L.I.}\fnms{L.I.}~\snm{Galtchouk}\corref{}\thanksref{A}\ead[label=e1]{leonid.galtchouk@math.unistra.fr}}
\and
\author[B]{\inits{S.M.}\fnms{S.M.}~\snm{Pergamenshchikov}\thanksref{B}\ead[label=e2]{Serge.Pergamenchtchikov@univ-rouen.fr}}
\address[A]{IRMA, Strasbourg University, 7 rue Rene Descartes,
67084, Strasbourg, France.\\ \printead{e1}}
\address[B]{Laboratoire de Math\'ematiques Raphael Salem,
Universit\'e de Rouen,
Avenue de l'Universit\'e, BP. 12,
F76801, Saint Etienne du Rouvray, Cedex France
and
Laboratory of Quantitative Finance,
National Research University -- Higher School of Economics,
Moscow, Russia.\\ \printead{e2}}
\end{aug}

\received{\smonth{3} \syear{2012}}
\revised{\smonth{12} \syear{2012}}

\begin{abstract}
A truncated sequential procedure
is constructed
for estimating the drift coefficient at a given state point based on
discrete data of ergodic
diffusion process.
A nonasymptotic upper bound is obtained for a pointwise absolute error risk.
The optimal
convergence rate and a sharp constant in the bounds are found for the
asymptotic pointwise
minimax risk. As a consequence, the efficiency is obtained of the
proposed sequential procedure.
\end{abstract}

\begin{keyword}
\kwd{discrete data}
\kwd{drift coefficient estimation}
\kwd{efficient procedure}
\kwd{ergodic diffusion process}
\kwd{minimax}
\kwd{nonparametric sequential
estimation}
\end{keyword}
\end{frontmatter}

\section{Introduction}\label{secIn}

In this paper, we consider
the following diffusion model:
%
\begin{equation}
\label{secIn.1}
\mathrm{d}y_{t}=S(y_{t})\,\mathrm{d}t +
\sigma(y_{t})\,\mathrm{d}W_{t},\qquad 0\le t\le T,
\end{equation}
where $(W_{t})_{t\ge0}$ is a scalar
standard Wiener process,
$S(\cdot)$ and $\sigma(\cdot)$ are unknown functions.
This model appears in a number of applied problems of stochastic
control, filtering,
spectral analysis,
identification of dynamic system, financial mathematics and
others (see \cite{Arato1982,Bensoussan1992,KaratzasShreve1998,Kutoyants2004,LiptserShiryaev1978}
and others for details).

The problem is to estimate the function $S(\cdot)$ at a point $x_{0}$
basing on the discrete time observations
%
\begin{equation}
\label{secIn.2}
(y_{t_{j}})_{1\le j\le N},\qquad  t_{j}=j\delta,
\end{equation}
where $N=[T/\delta]$ and the frequency $\delta=\delta_{T}\in
(0,1)$ is a
function of $T$ that will be specified later.

The estimation problem of the function $S$ was studied in a number of
papers in the case of complete
observations, that is when a continuous trajectory
$(y_{t})_{0\le t\le T}$ was observed. In the parametric case,
this problem
was considered apparently for the first time in the paper
\cite{AratoKolmogorovSinai1962} for diffusion model of the axis of the
equator precession.
In that paper, a nonasymptotic distribution of the maximum likelihood estimator
was found for a special Ornstein--Uhlenbeck process.

It should be noted that investigating nonasymptotic properties of
parametric estimators in the models like to \eqref{secIn.1} comes to
the analysis of nonlinear functionals of observations. At the most
cases, this analysis is unproductive in nonasymptotic setting. In
order to overcome the technique difficulties, the sequential analysis
methods were used in \cite{LiptserShiryaev1978} and \cite
{Novikov1972} for estimating a scalar parameter.
In \cite{KonevPergamenshchikov1981b}, these methods were extended to estimating
a multi-dimensional parameter as well. Moreover, in \cite{KonevPergamenshchikov1992}
truncated sequential procedures were developed that economizes the
observation time.

In \cite{GaltchoukPergamenshchikov2001} and \cite{GaltchoukPergamenshchikov2005},
a sequential approach was proposed for the pointwise nonparametric
estimation in the ergodic models \eqref{secIn.1}.
Later in \cite{GaltchoukPergamenshchikov2006b},
the efficiency was studied of the proposed sequential procedures.

A sufficiently complete survey one can find in \cite{Kutoyants2004} on
the nonparametric estimation in
the ergodic model \eqref{secIn.1}
when non sequential approaches are used.

In the cited papers, estimation problems were studied based on complete
observations $(y_{t})_{0\le t\le T}$.
In practice, usually one has at disposal discrete time observations
even for continuous time models.

A natural question arises about proprieties and the behavior of
estimates based on discrete
time observations for such models.
These problems were studied for several models. We cite some of them.

In \cite{Kessler1997}, asymptotically normal estimators are
constructed for both parameters $(\theta,\sigma)$, for a parametric
ergodic one-dimensional diffusion model observed at discrete times
$i\delta_n, 0\le i\le n$, with drift $b(x,\theta)$ and diffusion coefficient
$a(x,\sigma)$, when the observation frequency $\delta_n\to0$ and
$n\delta_n\to\infty$.
Estimating is based on the property of discrete observations to be
locally Gaussian. The author claims that asymptotic efficiency can be obtained
under additional condition \mbox{$n\delta_n^p\to0$}, where $p>1$.

Parametric estimations were studied in the papers \cite{Gobet2002,Yoshida1992}
and \cite{Jacod2006} for drift and diffusion coefficients in
multidimensional diffusion processes when the observation frequency
$\delta_n$ is as follows: $\delta_n\to0, n\delta_n\to\infty$.
In \cite{Gobet2002} the LAN property is proved in the ergodic case.
The proof is based on the transformation of the log-likelihood ratio by
the Malliavin calculus. Asymptotic normality is studied in \cite
{Yoshida1992} for the joint distribution of the maximum likelihood
estimator of parameters
in drift and diffusion coefficients.
In \cite{Jacod2006} the tightness of estimators is proved without
ergodicity or even recurrence assumptions.

Nonparametric estimation
setting for models of kind \eqref{secIn.1} was considered firstly
for estimating the unknown diffusion coefficient $\sigma^2(\cdot)$
based on discrete
time observations on a fixed interval $[0,T]$, when the observation
frequency goes to zero
(see, e.g., \cite{FlorensZmirou1993,GenonCatalotLaredoPicard1992,Hoffmann1999,Jacod2000} and the references therein).

Later, in
\cite{GobetHoffmannReiss2004} kernel estimates of drift and diffusion
coefficients were
studied for reflecting ergodic processes \eqref{secIn.1} taking the
values into the interval $[0,1]$ in the case of fixed
observation frequency; the asymptotics are taken as the sample size
goes to infinity. Minimax optimal convergence rates are found
for estimators of the drift and the diffusion coefficients. Upper and
lower bounds for $L_2$-risk are given as well.

So far as concerning the estimation in ergodic case, it should be noted
that a sequential
procedure was proposed
in \cite{Hoffmann1999} for nonparametric estimating the drift
coefficient of
the process \eqref{secIn.1} in the integral metric.
Some upper and lower asymptotic bounds were found for the $\mathbf
{L}_{p}$-risks.
Later, in the paper \cite{ComteGenenCatalotRozenholc2009} a nonasymptotic
oracle inequality was proved
for the drift coefficient estimation problem in a special empiric quadratic
risk based on discrete time observations. In the asymptotic setting,
when the observation frequency goes to zero and the length
of the observation time interval tends to infinity, the constructed
estimators reach the minimax optimal convergence rates.

This paper deals with the drift coefficient efficient nonparametric estimation
at a given state point based on discrete time observations \eqref{secIn.2}
in the absolute error risk. The unknown diffusion coefficient is a
nuisance parameter. We find the minimax optimal convergence rate and
we study the lower bound normalized by this convergence rate in the
case when the frequency $\delta_T\to0$ and $T\delta_T\to0$ as $T\to
\infty$.

Our approach is based on the sequential analysis developed in the papers
\cite{GaltchoukPergamenshchikov2001,GaltchoukPergamenshchikov2005},
and \cite{GaltchoukPergamenshchikov2006b} for the nonparametric
estimation. This approach
makes possible to replace the denominator by a constant in a sequential
Nadaraya--Watson estimator.

Let us recall that in the case of complete observations (i.e., when
a whole trajectory is
observed) the sequential estimate efficiency was proved by making use
of a uniform
concentration inequality (see \cite{GaltchoukPergamenshchikov2007}),
besides an indicator kernel estimator and a
weak H\"{o}lder space of functions $S$ were used.

As it turns out later in \cite{GaltchoukPergamenshchikov2011}, the efficient
kernel estimate in the above given sense provides to construct a
selection model adaptive
procedure that appears efficient in the quadratic $\mathbf{L}_{2}$-metric.

Therefore, in order to realize this program (i.e.,  from efficient
pointwise estimators to an efficient $\mathbf{L}_{2}$-estimator)
in the case of discrete time observations, one needs to obtain a
suitable concentration inequality, that is done in
\cite{GaltchoukPergamenshchikov2013}. It should be noted that in order
to obtain nonasymptotic
concentration inequality we make use of nonasymptotic bounds uniform
over functions
$S$ and $\sigma$ for the convergence rate in the ergodic theorem for
the process \eqref{secIn.1}. The latter result is proved in \cite{GaltchoukPergamenshchikov2014}
and it is based on a new approach using Lyapounov's functions and the
coupling method.

Further in this paper, we make use of the concentration inequality in
order to find the
explicit constant in the upper bound for weak H\"{o}lder's risk
normalized by the
optimal convergence rate and we prove that this upper bound is best
over all
possible estimators. It means the procedure is efficient.

The paper is organized as follows.
In Section~\ref{secFc}, we describe
the functional classes.
In Section~\ref{secSp} the sequential procedure is constructed.
In Section~\ref{secNa}, we obtain a nonasymptotic upper bound for the
absolute error pointwise risk
of the sequential procedure. In Section~\ref{secAe}, we show that the proposed
procedure is asymptotically efficient for the pointwise risk.
All proofs are given in Section~\ref{secPr}.
In the \hyperref[secApp]{Appendix}, we give all necessary technical results.

\section{Functional class}\label{secFc}

We consider the pointwise estimation problem
for the function $S(\cdot)$
at a fixed point $x_{0}\in\mathbb{R}$ for the model \eqref{secIn.1}
with unknown diffusion coefficient $\sigma$.
It is clear that to obtain a good estimate for the function $S(\cdot)$
at the point
$x_{0}$ it is necessary to impose some conditions
on the function $\vartheta=(S,\sigma)$
which provide that the observed process
$(y_{t})_{0\le t\le T}$ returns to any vicinity of the point
$x_{0}$
infinitely many times.

In this section, we describe a weak H\"older functional
class which guarantees the
ergodicity property for this model. First, for some $\mathbf{x}_{*}\ge|x_{0}|+1$, $M>0$ and $L>1$
we denote by $\Sigma_{L,M}$ the class of functions $S$ from
$\mathbf{C}^{1}(\mathbb{R})$ such that
\[
\sup_{|x|\le\mathbf{x}_{*}} \bigl(\bigl|S(x)\bigr|+\bigl|\dot{S}(x)\bigr| \bigr)\le M
\]
and
\[
-L\le\inf_{|x|\ge\mathbf{x}_{*}}\dot{S}(x)\le\sup_{|x|\ge\mathbf{x}_{*}}
\dot{S}(x) \le-L^{-1}.
\]

Moreover, for some fixed parameters
$0<\sigma_{\min}\le\sigma_{\max}<\infty$ we denote by
$\mathcal{V}$
the class of the functions $\sigma$ from $\mathbf{C}^{2}(\mathbb
{R})$ such that
%
\begin{equation}
\label{secFc.1}
\inf_{x\in\mathbb{R}} \bigl|\sigma(x)\bigr| \ge\sigma_{\min}
\quad\mbox{and}\quad\sup_{x\in\mathbb{R}} \max \bigl(\bigl|\sigma(x)\bigr|, \bigl|\dot{\sigma}(x)\bigr|,
\bigl|\ddot{\sigma }(x)\bigr| \bigr) \le\sigma_{\max}.
\end{equation}

In this paper, we make use of the weak H\"{o}lder functions introduced in
\cite{GaltchoukPergamenshchikov2006a}.

\begin{definition}\label{De.secFc.1}
We say that a function $S$ satisfies the
weak H\"{o}lder condition at the point $x_{0}\in\mathbb{R}$
with the parameters $h,\varepsilon>0$
and exponent $\beta=1+\alpha$, $\alpha\in(0,1)$, if
$S\in\mathbf{C}^1(\mathbb{R})$ and its derivative satisfies the following
inequality
%
\begin{equation}
\label{secFc.2}
\biggl\llvert \int^1_{-1} z\int
^1_0 \bigl(\dot{S}(x_{0}+uzh)-
\dot{S}(x_{0})\bigr) \,\mathrm{d}u \,\mathrm{d}z\biggr\rrvert \le\varepsilon
h^{\alpha}.
\end{equation}
We will denote the set of all such functions by
$\mathcal{H}^{w}_{x_{0}}(\varepsilon,\beta,h)$.
\end{definition}

Note that
%
\begin{equation}
\label{secFc.2-1}
\int^1_{-1} z\int
^1_0 \bigl(\dot{S}(x_{0}+uzh)-
\dot{S}(x_{0})\bigr) \,\mathrm{d}u \,\mathrm{d}z=\Omega_{x_{0},h}(S),
\end{equation}
where
$\Omega_{x_{0},h}(S)=\int^{1}_{-1}
 (
S(x_{0}+hz)
-
S(x_{0})
 )\,
\mathrm{d}z$.
Therefore, the condition \eqref{secFc.2} for the functions from
$\mathcal{H}^{w}_{x_{0}}(\varepsilon,\beta,h)$ is equivalent to
the following one
%
\begin{equation}
\label{secFc.3}
\sup_{S\in\mathcal{H}^{w}_{x_{0}}(\varepsilon,\beta,h)} \bigl\llvert \Omega_{x_{0},h}(S)
\bigr\rrvert \le \varepsilon h^{\beta}.
\end{equation}
Let us
denote by $\mathcal{H}^{w}_{x_{0},M}(\varepsilon,\beta,h)$
the set of all functions $D$ from $\mathcal{H}^{w}_{x_{0}}(\varepsilon,\beta,h)$
such that $\sup_{x\in\mathbb{R}}(|D(x)|+|\dot{D}(x)|)\le M/2$ and
$D(x)=0$ for
$|x|\ge x_{*}$.

Let $S_{0}$ be a function from
$\Sigma_{L,M/2}$ such that
%
\begin{equation}
\label{secFc.4}
\lim_{h\to0} h^{-\beta}
\Omega_{x_{0},h}(S_{0})=0.
\end{equation}
We define a vicinity $\mathcal{U}_{M}(x_{0},\beta)$ of the
function $S_{0}$ as follows
%
\begin{equation}
\label{secFc.5}
\mathcal{U}_{M}(x_{0},\beta)=
S_{0}+ \mathcal{H}^{w}_{x_{0},M}(\varepsilon,\beta,h),
\end{equation}
where $h=T^{-1/(2\beta+1)}$
and
%
\begin{equation}
\label{secFc.6}
\varepsilon=\varepsilon_{T}= \frac{1}{(\ln T)^{1+\gamma}}
\end{equation}
for some $0<\gamma<1$.
Obviously that $\mathcal{U}_{M}(x_{0},\beta)\subset\Sigma
_{L,M}$.
Finally, we set
%
\begin{equation}
\label{secFc.9}
\Theta_{\beta}= \mathcal{U}_{M}(x_{0},
\beta)\times\mathcal{V}.
\end{equation}
It should be noted that, for any $\vartheta\in\Theta_{\beta}$,
there exists an invariant density which is defined as
%
\begin{equation}
\label{secFc.10}
q_{\vartheta}(x)= \biggl(\int_{\mathbb{R}}
\sigma^{-2}(z) \mathrm{e}^{\widetilde{S}(z)}\,\mathrm{d}z \biggr)^{-1}
\sigma^{-2}(x) \mathrm{e}^{\widetilde{S}(x)},
\end{equation}
where $\widetilde{S}(x)=2\int_{0}^{x} \sigma^{-2}(v)S(v)\,\mathrm{d}v$
(see, e.g., \cite{GihmanSkorohod1968}, Chapter~4.18, Theorem~2).
It is easy to see that
this density is uniformly bounded in the class \eqref{secFc.9}, that is,
%
\begin{equation}
\label{secFc.11}
q^{*} = \sup_{x\in\mathbb{R}} \sup
_{\vartheta\in\Theta
_{\beta}} q_{\vartheta}(x) < +\infty
\end{equation}
and bounded away from zero on the interval $[x_{0}-1,x_{0}+1]$,
that is,
%
\begin{equation}
\label{secFc.12}
q_{*} = \inf_{x_{0}-1\le x\le x_{0}+1}\ \inf
_{\vartheta\in\Theta_{\beta}} q_{\vartheta}(x) > 0.
\end{equation}
For any $\mathbb{R}\to\mathbb{R}$ function $f$ from $\mathbf{L}_{1}(\mathbb{R})$, we set
%
\begin{equation}
\label{secFc.13}
\mathbf{m}_{\vartheta}(f)=\int_{\mathbb{R}} f(x)
q_{\vartheta}(x) \,\mathrm{d}x.
\end{equation}

Assume that the frequency $\delta$ in the observations
\eqref{secIn.2}
is of the following asymptotic form (as $T\to\infty$)
%
\begin{equation}
\label{secFc.14}
\delta=\delta_{T}=\mathrm{O} \biggl( \frac{\varepsilon_{T}}{T}
\biggr),
\end{equation}
where the function $\varepsilon_{T}$
is introduced in \eqref{secFc.6}.

Now, for any estimate (i.e., any
$(y_{t})_{0\le t\le T}$
measurable function) $\tilde{S}_{T}(x_{0})$ of $S(x_{0})$,
we define the pointwise risk as follows
%
\begin{equation}
\label{secAe.1}
\mathcal{R}_{\vartheta}(\widetilde{S}_{T})=
\mathbf{E}_{\vartheta} \bigl|\widetilde{S}_{T}(x_{0})-S(x_{0})\bigr|.
\end{equation}

\begin{remark}\label{Re.secFc.1}
It should be noted that in the definition of the weak H\"{o}lder class
$\mathcal{H}^{w}_{x_{0}}(\varepsilon,\beta,h)$ in
\eqref{secFc.2} the weak H\"{o}lder norm $\varepsilon$ is chosen such
that it goes to zero as $T\to\infty$
(see \eqref{secFc.6}) in contrast with the initial definition of this
class given in
\cite{GaltchoukPergamenshchikov2006b}, where the norm was fixed. But
there an
additional limit passage with $\varepsilon\to0$ was done in the theorem about
upper bound (see, Theorem~5.2 in \cite{GaltchoukPergamenshchikov2006b}).
Therefore, in this paper we choose $\varepsilon\to0$ in the definition
of the weak H\"{o}lder class
instead of the additional limit passage with $\varepsilon\to0$ in the
theorem on the upper bound.
Technically it is nearly the same, but the sense of the upper bound is
clearer without
the additional limit with $\varepsilon\to0$.
\end{remark}

\section{Sequential procedure}\label{secSp}

In order to construct an efficient pointwise estimator of $S$,
we begin with estimating the ergodic
density $q=q_{\vartheta}$ at the point $x_{0}$
from first $N_{0}$ observations. We choose
%
\begin{equation}
\label{secSp.0}
N_{0}=N^{\gamma_{0}}\quad \mbox{and}\quad 2/3<
\gamma_{0}<1.
\end{equation}
We will make use of the following kernel estimator
%
\begin{equation}
\label{secSp.1}
\widehat{q}_{T}(x_{0})=\frac{1}{2(N_{0}-1)\varsigma}
\sum^{N_{0}-1}_{j=0} Q \biggl(\frac{y_{t_{j}}-x_{0}}{\varsigma}
\biggr),
\end{equation}
where $Q(y)=\mathbf{1}_{ (|y|\le1 )}$ and $\varsigma
=\varsigma_{T}$
is a function of $T$
such that
\[
\varsigma_{T}=\mathrm{o}\bigl(T^{-\gamma_{0}/2}\bigr)\qquad \mbox{as } T\to
\infty.
\]

For $T\ge3$, we set
%
\begin{equation}
\label{secSp.1-1}
\widetilde{q}_{T}(x_{0})=
\cases{(\upsilon_{T})^{1/2},&\quad $\mbox{if }\widehat{q}_{T}(x_{0})< (\upsilon_{T})^{1/2}$;
\vspace*{3pt}\cr
\widehat{q}_{T}(x_{0}),&\quad$\mbox{if } (\upsilon_{T})^{1/2}
\le\widehat{q}_{T}(x_{0})\le (\upsilon_{T})^{-1/2}$;
\vspace*{3pt}\cr
(\upsilon_{T})^{-1/2},&\quad$\mbox{if } \widehat{q}_{T}(x_{0})>
(\upsilon_{T})^{-1/2}$,}
\end{equation}
where
\[
\upsilon_{T}=\frac{1}{(\ln T)^{a_{0}}}\quad \mbox{and}\quad a_{0}=
\frac{\sqrt{\gamma+1}-1}{10}.
\]
The properties of the estimates $\widehat{q}_{T}(x_{0})$ and
$\widetilde{q}_{T}(x_{0})$
are studied in the \hyperref[secApp]{Appendix}.

Let us define the following stopping time
%
\begin{equation}
\label{secSp.2}
\varpi=\varpi_{T}=\inf \Biggl\{j\ge N_{0} \dvt
\sum^{j}_{i=N_{0}} \phi_{i} \ge
H_{T} \Biggr\},
\end{equation}
where $H_{T}$ is a threshold,
$
\phi_{i}=\chi_{h,x_{0}}(y_{t_{i-1}})\mathbf
{1}_{\{i\le N\}}
+
\mathbf{1}_{\{i> N\}}$,
$\chi_{h,x_{0}}(y)=Q ((y-x_{0})/h )$
and $h$ is a positive bandwidth. We put
$\varpi=\infty$ if the set $ \{\cdot \}$ is empty.
Obviously, that in the our case
$\varpi<\infty$ a.s. since $\sum_{i\ge N_{0}} \phi_{i}=+\infty$ a.s.

Now we have to choose the threshold $H_{T}$. Note that in order to
construct an
efficient estimator one should use all, that is, $N$ observations.
Therefore, the threshold $H_{T}$ should provide the asymptotic
relations $\varpi_{T}\approx N$
and
\[
\sum^{N}_{i= N_{0}} \phi_{i} = \sum
^{N}_{i= N_{0}} \chi_{h,x_{0}}(y_{t_{i-1}})
\qquad\mbox{as }T\to\infty.
\]
In order to obtain these relations, note that, due to the ergodic theorem,
\[
\sum^{N}_{i= N_{0}} \chi_{h,x_{0}}(y_{t_{i-1}})
\approx 2 h (N-N_{0}) q_{\vartheta}(x_{0}).
\]
Hence, replacing in the right-hand side term the ergodic density with its
corrected estimate yields the following definition of the threshold
%
\begin{equation}
\label{secSp.3}
H=H_{T}=h(N-N_{0}) \bigl(2
\widetilde{q}_{T}(x_{0})-\upsilon _{T}\bigr).
\end{equation}
Note that in \cite{GaltchoukPergamenshchikov2001} it has been shown
that the such form of the threshold $H_{T}$ provides the optimal
convergence rate.
It is clear that
\[
\varpi\le N+H_{T}<N+h (N-N_{0})/\sqrt{
\upsilon_{T}},
\]
that is,
the stopping time $\varpi$ is bounded.
Now on the set $\Gamma_{T}= \{\varpi\le N \}$ we
define the correction
coefficient
$\varkappa=\varkappa_{T}$ as
\[
\varkappa_{T}=\frac
{H_{T}-\sum^{\varpi-1}_{j=N_{0}} \chi_{h,x_{0}}(y_{t_{j-1}})}{
\chi_{h,x_{0}}(y_{t_{\varpi-1}})},
\]
that is, on the set $\Gamma_{T}$
\[
\sum^{\varpi-1}_{j=N_{0}} \chi_{h,x_{0}}(y_{t_{j-1}})+
\varkappa\chi_{h,x_{0}}(y_{t_{\varpi-1}})=H_{T}.
\]
Moreover, on the $\Gamma^{c}_{T}$ we set $\varkappa_{T}=1$.
Using this definition,
we introduce the weight sequence
%
\begin{equation}
\label{secSp.4}
\widetilde{\varkappa}_{j}=\mathbf{1}_{\{j<\varpi\}}+
\varkappa \mathbf{1}_{\{j=\varpi\}},\qquad j\ge1.
\end{equation}
One can check directly that, for any $j\ge1$, the coefficients
$\widetilde{\varkappa}_{j}$ are
$\mathcal{F}_{t_{j-1}}$ measurable, where
$\mathcal{F}_{t_{j}}=\sigma (y_{t_{k}}, 0\le
k\le j )$.
Now we define the sequential estimator for $S(x_{0})$ as
%
\begin{equation}
\label{secSp.5}
S^{*}_{h,T}(x_{0})=
\frac{1}{\delta H_{T}} \Biggl( \sum^{N}_{j=N_{0}}
\widetilde{\varkappa}_{j} \chi_{h,x_{0}}(y_{t_{j-1}}) \Delta
y_{t_{j}} \Biggr) \mathbf{1}_{\Gamma_{T}}.
\end{equation}
In the next section, we study nonasymptotic
properties of this procedure.

\begin{remark}\label{Re.secSp.1}
Note that the correction coefficient of type \eqref{secSp.4} was used first
in the paper \cite{BorisovKonev1977} in order to construct an unbiased
estimator
of a scalar parameter in autoregressive processes AR(1). Here, we make
use of the
same idea for a nonparametric procedure.
\end{remark}

\begin{remark}\label{Re.secSp.2}
In fact, our procedure uses only the observations belonging to the
interval $[x_{0}-h, x_{0}+h]$,
that results in the sample size asymptotically equals to $2Nh\widehat
{q}_{T}(x_{0})$.
This is related with the choice of the estimator kernel that is an
indicator function. It is easy to verify
(see \cite{GaltchoukPergamenshchikov2005}) that this kernel minimizes
the variance of stochastic term
in the kernel estimator. Ultimately, the last result provides
efficiency of the procedure.
\end{remark}

\section{Nonasymptotic estimation}\label{secNa}

In this section, an upper bound for the absolute error risk will be
given in the case, when $S\in\Sigma_{L,M}, \sigma$
is differentiable and $\sigma_{\min}\le|\sigma(x)|\le\sigma_{\max}$ for any $x$. We will denote this case as
$\vartheta\in\Sigma_{L,M}\times[\sigma_{\min},\sigma_{\max}]$.

As we will see later in studying the estimator \eqref{secSp.5}, the
approximation term
plays the crucial role. In the our case, this term is of the following form
%
\begin{equation}
\label{secUp.1}
\Upsilon_{1,T}=\frac{1}{\delta H_{T}} \sum
^{N}_{j=N_{0}} \widetilde{\varkappa}_{j}
\chi_{h,x_{0}}(y_{t_{j-1}}) \varrho_{j},
\end{equation}
where
$\varrho_{j}=
\int^{t_{j}}_{t_{j-1}}(S(y_{u})-S(y_{t_{j-1}}))\,\mathrm{d}u$.
One can show the following result.

\begin{proposition}\label{Pr.secUp.1}
For any $T\ge3$,
%
\begin{equation}
\label{secUp.2}
\sup_{\vartheta\in\Sigma_{L,M}\times[0,\sigma_{\max}]} \mathbf{E}_{\vartheta}
\Upsilon^{2}_{1,T} \le \widetilde{L}^{2}L_{1}
\delta,
\end{equation}
where $\widetilde{L}=\max(L,M)$
and
$L_{1}=2 (\sigma^{2}_{\max}+2\delta(M^{2}+L^{3}D_{*}+L^{2}
x^{2}_{*}) )$.
\end{proposition}

Further, we set
%
\begin{equation}
\label{secUp.3} \Upsilon_{2,T}=\frac{1}{\delta H_{T}} \sum
^{N}_{j=N_{0}} \widetilde{\varkappa}_{j}
\chi_{h,x_{0}}(y_{t_{j-1}}) \varrho^{*}_{j},
\end{equation}
where
$\varrho^{*}_{j}=
\int^{t_{j}}_{t_{j-1}}(\sigma(y_{u})-\sigma(y_{t_{j-1}}))\,\mathrm{d}W_{u}$.

\begin{proposition}\label{Pr.secUp.2}
For any $T\ge3$ for which $0<\delta\le1$,
one has
%
\begin{equation}
\label{secUp.4}
\sup_{\vartheta\in\Sigma_{L,M}\times[0,\sigma_{\max}]} \mathbf{E}_{\vartheta} (
\Upsilon_{2,T} )^{2} \le \frac{\sigma^{2}_{\max}
L_{1}}{h(N-N_{0})\sqrt{\upsilon_{T}}}.
\end{equation}
\end{proposition}

Proofs of Propositions~\ref{Pr.secUp.1} and \ref{Pr.secUp.2}
are given in the \hyperref[secApp]{Appendix}.

Now we introduce the approximative term, that is,
%
\begin{equation}
\label{secUp.4-0}
B_{T}= \frac{1}{H_{T}} \sum
^{N}_{j=N_{0}} \widetilde{\varkappa}_{j}
f_{h}(y_{t_{j-1}})
\end{equation}
with $f_{h}(y)=\chi_{h,x_{0}}(y)(S(y)-S(x_{0}))$.
Taking into account this formula, we can represent
the error of estimator \eqref{secSp.5}
on the set $\Gamma_{T}$ as
%
\begin{equation}
\label{secUp.4-0-1}
S^{*}_{h,T}(x_{0})-S(x_{0})=
\Upsilon_{1,T} + B_{T}+ \mathbf{M}_{T},
\end{equation}
where
\[
\mathbf{M}_{T}=\frac{1}{\delta H_{T}}\sum^{N}_{j=N_{0}}
\widetilde{\varkappa}_{j} \chi_{h,x_{0}}(y_{t_{j-1}})
\eta_{j}
\]
with $\eta_{j}=\int^{t_{j}}_{t_{j-1}}\sigma(y_{u}) \,\mathrm{d}W_{u}$.
Obviously, for any function $S$ from $\Sigma_{L,M}$, the term
$B_{T}$ can be bounded
as
%
\begin{equation}
\label{secNa.0}
|B_{T}|\le h \max_{|x-x_{0}|\le h}\bigl|\dot{S}(x)\bigr|
\le M h.
\end{equation}

\begin{proposition}\label{Pr.secUp.3}
For any $T\ge3$,
one has
%
\begin{equation}
\label{secUp.4-2}
\sup_{\vartheta\in\Sigma_{L,M}\times[0,\sigma_{\max}]} \mathbf{E}_{\vartheta}
\mathbf{M}^{2}_{T} \le \frac{\sigma^2_{\max}
}{\delta h(N-N_{0})\sqrt{\upsilon_{T}}}.
\end{equation}
\end{proposition}

Hence, we obtain the following upper bound.

\begin{theorem}\label{Th.secNa.1}
For any $h>0$ and $T\ge3$ for which $0<\delta\le1$,
one has
%
\begin{equation}
\label{secNa.1}
\sup_{\vartheta\in\Sigma_{L,M}\times[\sigma_{\min},
\sigma_{\max}]} \mathbf{E}_{\vartheta} \bigl|
S^{*}_{h,T}(x_{0})-S(x_{0}) \bigr| \le
U^{*}(\delta,h,T) + M \Pi^{*}_{T},
\end{equation}
where
\[
U^{*}(\delta,h,T)=\widetilde{L}\sqrt{\delta L_{1}} + M h +
\frac{\sigma_{\max}}{\sqrt{\delta h (N-N_{0})}
\upsilon^{1/4}_{T}}
\]
and
\[
\Pi^{*}_{T}= \sup_{\vartheta\in\Sigma_{L,M}\times[\sigma_{\min},
\sigma_{\max}]}
\mathbf{P}_{\vartheta} \bigl( \Gamma^{c}_{T} \bigr).
\]
\end{theorem}

Let us study now the last term in \eqref{secNa.1}.

\begin{proposition}\label{Pr.secNa.1}
Assume that the parameter $\delta$ is of the asymptotic form
\eqref{secFc.14} and $h\ge T^{-1/2}$. Then, for any $a>0$,
%
\begin{equation}
\label{secSp.5-1}
\lim_{T\to\infty} T^{a}
\Pi^{*}_{T} =0.
\end{equation}
\end{proposition}

Proof of this proposition is given in the \hyperref[secApp]{Appendix}.

\begin{remark}\label{Re.secNa.1}
It should be noted that the main destination of the bound \eqref
{secNa.1} is
to obtain a sharp oracle inequality for a model selection procedure for
the process
\eqref{secIn.1} observed at discrete times.
Recall (see, \cite{GaltchoukPergamenshchikov2011}), that an efficient
model selection procedure
for diffusion processes is based on estimators admitting on some set
$\Gamma_{T}$ a nonasymptotic
representation of kind \eqref{secUp.4-0-1} satisfying the conditions
\eqref{secNa.0} and
\eqref{secUp.4-2} at estimation state points $(x_{j})$. In
addition, the condition \eqref{secSp.5-1}
must hold true for the set $\Gamma_{T}$. Moreover, the stochastic
terms of kernel estimators,
$\mathbf{M}_{T}=\mathbf{M}_{T}(x_{j})$, must be
independent random variables for different points $x_{j}$.
The last condition is provided by sequential approach since, for
sequential kernel estimator, the term
$\mathbf{M}_{T}$ is a Gaussian random variable on the set $\Gamma
_{T}$.
\end{remark}

\section{Asymptotic efficiency}\label{secAe}

First of all, we study a lower bound for the risk \eqref{secAe.1}. To
this end,
we set
%
\begin{equation}
\label{secAe.2}
\varsigma^{*}_{\vartheta} (x_{0})=
\frac{2q_{\vartheta}(x_{0})}{\sigma^2(x_{0})}.
\end{equation}
This parameter provides a sharp asymptotic lower bound for the pointwise
risk normalized by the minimax rate
$\varphi_{T}=T^{\beta/(2\beta+1)}$.

\begin{theorem}\label{Th.secLo.1}
The risk defined in \eqref{secAe.1} admits the following lower bound
%
\begin{equation}
\label{secAe.3}
\mathop{\underline{\lim}}_{T\to\infty}\varphi_{T}\inf
_{\widetilde
{S}_{T}} \sup_{\vartheta\in\Theta_{\beta}} \sqrt{
\varsigma^{*}_{\vartheta}(x_{0})} \mathcal{R}_{\vartheta}(
\widetilde{S}_{T}) \ge \mathbf{E}|\xi|,
\end{equation}
where infimum is taken over all possible estimators $\widetilde{S}_{T}$,
$\xi$ is a $(0,1)$ Gaussian random variable.
\end{theorem}

\begin{theorem}\label{Th.secUp.1}
The kernel estimator $S^{*}_{h,T}$ defined in \eqref{secSp.5}
with
$h=T^{-1/(2\beta+1)}$
satisfies the following asymptotic inequality
%
\begin{equation}
\label{secAe.3-0}
\mathop{\overline{\lim}}_{T\to\infty} \varphi_{T} \sup
_{\vartheta\in\Theta_{\beta}} \sqrt{\varsigma^{*}_{\vartheta}(x_{0})}
\mathcal{R}_{\vartheta}\bigl(S^{*}_{h,T}\bigr) \le
\mathbf{E}|\xi|,
\end{equation}
where $\xi$ is a $(0,1)$ Gaussian random variable.
\end{theorem}
%

%
\begin{remark}\label{Re.secUp.0}
It should be noted that one cannot use the bound \eqref{secNa.1} in
order to obtain the
sharp asymptotic upper bound \eqref{secAe.3-0} because the upper
bound \eqref{secNa.1}
is obtained for a wider function class,
that is, for $\vartheta\in\Sigma_{L,M}\times[\sigma_{\min}, \sigma_{\max}]$
and hence, it is not the best.
\end{remark}

Notice that the Theorems~\ref{Th.secLo.1}
and~\ref{Th.secUp.1} imply the following
efficiency property.

\begin{theorem}\label{Re.secEff.1}
The sequential procedure
\eqref{secSp.5}
with
$h=T^{-1/(2\beta+1)}$
is
asymptotically efficient in the following sense:
\[
\mathop{\overline{\lim}}_{T\to\infty} \varphi_{T} \sup_{\vartheta\in\Theta_{\beta}}
\sqrt{\varsigma^{*}_{\vartheta}(x_{0})}
\mathcal{R}_{\vartheta}\bigl(S^{*}_{h,T}\bigr) =
\mathop{\underline{\lim}}_{T\to\infty}\varphi_{T}\inf_{\widetilde
{S}_{T}}
\sup_{\vartheta\in\Theta_{\beta}} \sqrt{\varsigma^{*}_{\vartheta}(x_{0})}
\mathcal{R}_{\vartheta}(\widetilde{S}_{T}),
\]
where infimum is taken over all possible estimators $\widetilde{S}_{T}$.
\end{theorem}

\begin{remark}\label{Re.secUp.1}
The constant \eqref{secAe.2} provides the sharp asymptotic lower
bound for the minimax
pointwise risk. The calculation of this
constant is possible by making use of the weak H\"{o}lder class. This
functional class was introduced in \cite{GaltchoukPergamenshchikov2006a} for regression models.
For the first time, the constant \eqref{secAe.2} was obtained in the paper
\cite{GaltchoukPergamenshchikov2006b} at the pointwise estimation
problem of the drift based on continuous time observations of the
process \eqref{secIn.1} with the unit diffusion. Later this constant
was used in the paper \cite{GaltchoukPergamenshchikov2011} to obtain
the Pinsker constant for a quadratic risk in the adaptive estimation
problem of the drift in the model \eqref{secIn.1} based on continuous
time observations.
\end{remark}

\begin{remark}\label{Re.secUp.2}
Note also that in this paper the efficient procedure is constructed when
the regularity is known of the function to be estimated.
In the case of unknown regularity, we shall use an approach based on
the model selection similarly to that in the
paper \cite{GaltchoukPergamenshchikov2011} which deals with continuous
time observations.
The announced result will be published in the next paper which is in
the work.
\end{remark}

\begin{remark}\label{Re.secUp.3}
In the paper, we studied only the case of H\"{o}lderian smoothness
$1+\alpha$ with $\alpha\in(0,1)$.
Efficiency is provided by the indicator kernel which minimizes the
asymptotic variance of the stochastic
term (see, e.g., \cite{GaltchoukPergamenshchikov2005}). If in
the pointwise estimation problem
the unknown function possesses a greater smoothness then, as known, one
needs to make use of a kernel that
should be orthogonal to all polynomials of orders less than the
integral part of the smoothness order.
It is clear that the kernel estimator is not efficient. Therefore, one
needs to uses an other estimator, in
particular, a local polynomial estimator but again with an indicator kernel.
This is a subject of our future investigation in the pointwise setting
for diffusion processes. It will
based on ideas and results of $1+\alpha$-smoothness case but, due to
extreme complication, the new case
cannot be considered in this paper.
\end{remark}

\section{Proofs}\label{secPr}

\subsection{Lower bound}\label{secLo}

In this section, the Theorem~\ref{Th.secLo.1} will be proved.
Let us introduce the model \eqref{secIn.1} with $\sigma=1$, that is,
%
\begin{equation}
\label{secLo.1}
\mathrm{d}y_{t}=S(y_{t}) \,\mathrm{d}t +
\mathrm{d}W_{t}.
\end{equation}
Now we define the risk corresponding to this model as follows
%
\begin{equation}
\label{secLo.2} \mathcal{R}^{*}_{S}(\widetilde{S}_{T})=
\mathbf{E}_{S} \bigl|\widetilde{S}_{T}(x_{0})-S(x_{0})\bigr|,
\end{equation}
where $\mathbf{E}_{S}$ denotes the expectation with respect to the
distribution $\mathbf{P}_{S}$
of the process \eqref{secLo.1} in the space of continuous functions
$\mathbf{C}[0,T]$. It is clear that
%
\begin{equation}
\label{secLo.3}
\sup_{\vartheta\in\Theta_{\beta}} \sqrt{\varsigma^{*}_{\vartheta}(x_{0})}
\mathcal{R}_{\vartheta}(\widetilde{S}_{T}) \ge \sup
_{S\in\mathcal{U}_{M}(x_{0},\beta)} \sqrt{2 q_{S}(x_{0})}
\mathcal{R}^{*}_{S}(\widetilde{S}_{T}),
\end{equation}
where $q_{S}$ is the invariant density for the process \eqref{secLo.1} which equals
to $q_{\vartheta}$ with $\sigma=1$. Let now $g$ be a continuously
differentiable
probability density on the interval $[-1,1]$. Then, for any $u\in
\mathbb{R}$ and
$0<\nu<1/4$, we set
\[
S_{u,\nu}(x)=S_{0}(x)+\frac{u}{\varphi_{T}} V_{\nu}
\biggl(\frac{x-x_{0}}{h} \biggr),
\]
where $h=T^{-1/(2\beta+1)}$ and
\[
V_{\nu}(x)=\frac{1}{\nu}\int^{\infty}_{-\infty}
( \mathbf{1}_{ (|u|\le1-2\nu )}+ 2\mathbf{1}_{ (1-2\nu\le|u|\le1-\nu )} ) g \biggl(
\frac{u-x}{\nu} \biggr) \,\mathrm{d}u.
\]
It is easy to see directly that, for any $0<\nu<1/4$,
\[
V_{\nu}(0)=1 \quad\mbox{and}\quad \int^1_{-1}V_{\nu}(x)
\, \mathrm{d}x=2.
\]
Therefore, denoting $D_{u}(x)=S_{u,\nu}(x)-S_{0}(x)$, we
obtain, for any
$u\in\mathbb{R}$,
\[
\Omega_{x_{0},h}(D_{u})=\int^{1}_{-1}
\bigl(D_{u}(x_{0}+hz)-D_{u}(x_{0})
\bigr) \,\mathrm{d}z=0.
\]
Moreover, note that, for any fixed $b>0$,
\[
\sup_{|u|\le b} \bigl|\dot{D}_{u}(x)\bigr|= \sup
_{|u|\le b} \biggl( |u| \varphi_{T}^{-1}h^{-1}
\biggl|\dot{V}_{\nu} \biggl(\frac{x-x_{0}}{h} \biggr)\biggr| \biggr) \le b
T^{-\alpha/(2\beta+1)}\nu^{-2}\dot{g}^{*},
\]
where $\dot{g}^*=\sup_x|\dot{g}(x)|$.
Therefore, $\sup_{x\in\mathbb{R}}\sup_{|u|\le b}|\dot
{D}_{u}(x)|\le M/2$
for sufficiently large $T$ and, in view of the equality \eqref{secFc.2-1},
the functions
$(S_{u,\nu})_{|u|\le b}$ belong to the class $\mathcal{U}_{M}(x_{0},\beta)$
for sufficiently large $T$. It implies that, for any $b>0$ and for
sufficiently large
$T$, we can estimate from below
the right-hand term in the inequality \eqref{secLo.3} as
\begin{eqnarray*}
\sup_{S\in\mathcal{U}_{M}(x_{0},\beta)} \sqrt{2 q_{S}(x_{0})}
\mathcal{R}^{*}_{S}(\widetilde{S}_{T}) &\ge&  \sup
_{|u|\le b} \sqrt{2 q_{S_{u,\nu}}(x_{0})}
\mathcal{R}^{*}_{S_{u,\nu}}(\widetilde{S}_{T})
\\
&=& \sup_{|u|\le b} \sqrt{2 q_{S_{0}}(x_{0})}
\mathcal{R}^{*}_{S_{u,\nu}}(\widetilde{S}_{T})+Q_{T},
\end{eqnarray*}
where $Q_{T}=\sup_{|u|\le b}
\sqrt{2 q_{S_{u,\nu}}(x_{0})}
\mathcal{R}^{*}_{S_{u,\nu}}(\widetilde{S}_{T})
-\sup_{|u|\le b}
\sqrt{2 q_{S_{0}}(x_{0})}
\mathcal{R}^{*}_{S_{u,\nu}}(\widetilde{S}_{T})$.

It is easy to see that
\begin{eqnarray*}
|Q_{T}| &\le &
\sup_{|u|\le b} \bigl|\sqrt{2 q_{S_{u,\nu}}(x_{0})} -
\sqrt{2 q_{S_{0}}(x_{0})}\bigr| \mathcal{R}^{*}_{S_{u,\nu}}(
\widetilde{S}_{T}) \\
&\le & \frac{1}{\sqrt{2q_*}}\sup_{|u|\le b} \bigl|
q_{S_{u,\nu}}(x_{0}) - q_{S_{0}}(x_{0})\bigr|
\mathcal{R}^{*}_{S_{u,\nu}}(\widetilde{S}_{T}).
\end{eqnarray*}
Taking into account here that
\[
\lim_{T\to\infty} \sup_{|u|\le b} \bigl\llvert
q_{S_{u,\nu}}(x_{0}) - q_{S_{0}}(x_{0}) \bigr\rrvert = 0,
\]
we obtain the inequality \eqref{secAe.3}
by making use of the Theorem~4.1 from
\cite{GaltchoukPergamenshchikov2006b}.
Thus, we obtain the Theorem~\ref{Th.secLo.1}.

\subsection{Upper bound}\label{secUp}

We begin with stating the following result for the term \eqref{secUp.4-0}.

%
\begin{proposition}\label{Pr.secUp.4}
The function $B_{T}$ defined in \eqref{secUp.4-0} satisfies the
following asymptotic
property
%
\begin{equation}
\label{secUp.8} \limsup_{T\to\infty} \varphi_{T} \sup
_{\vartheta\in\Theta_{\beta}} \mathbf{E}_{\vartheta} |B_{T}|=0.
\end{equation}
\end{proposition}

The result is proved in the \hyperref[secApp]{Appendix}.

Now we prove Theorem~\ref{Th.secUp.1}. To this end, we set
\[
\widetilde{\phi}(u)= \sum_{j=N_{0}}^{+\infty}
\phi_{j} \mathbf{1}_{ \{t_{j-1}< u\le t_{j} \}},
\]
where the random variables $(\phi_{i})_{i\ge1}$ are defined in
\eqref{secSp.2}. Using this function, we introduce the stopping time
\[
\tau=\tau_{T}=\inf \biggl\{t\ge T_{0} \dvt \int
^{t}_{T_{0}} \widetilde{\phi}(u) \,\mathrm{d}u\ge\delta
H_{T} \biggr\},
\]
where $T_{0}=t_{N_{0}}=\delta N_{0}$.
As usually, we put $\tau=\infty$ if the set
$ \{\cdot \}$ is empty. Obviously that
\[
\tau\le T+\delta H_{T}\le T+\delta h(N-N_{0})/\sqrt{
\upsilon _{T}}.
\]
Due to the equality $\int^{\infty}_{T_{0}}
 \widetilde{\phi}(u)\, \mathrm{d}u=\infty$,
one obtains immediately that the random
variable
%
\begin{equation}
\label{secUp.6} \xi_{T}=\frac{1}{\sqrt{\delta H_{T}}} \int^{\tau}_{T_{0}}
\widetilde{\phi}(u) \,\mathrm{d}W_{u}
\end{equation}
is Gaussian $\mathcal{N}(0,1)$ (see, e.g., \cite
{LiptserShiryaev1978}, Chapter~17).
Now, using this property,
we can rewrite the deviation
\eqref{secUp.4-0-1}
on set $\Gamma_{T}$ as
%
\begin{equation}
\label{secUp.7}
S^{*}_{h,T}(x_{0})-S(x_{0})=
B^{*}_{T} + \mathbf{M}^{(1)}_{T} +
\sigma(x_{0}) \mathbf{M}^{(2)}_{T} +
\frac{\sigma(x_{0})}{\sqrt{\delta H_{T}}} \xi_{T},
\end{equation}
where $B^{*}_{T}=\Upsilon_{1,T}
+\Upsilon_{2,T}+B_{T}$,
\[
\mathbf{M}^{(1)}_{T} = \frac{1}{\delta H_{T}} \sum
^{N}_{j=N_{0}} \widetilde{\varkappa}_{j}
\chi_{h,x_{0}}(y_{t_{j-1}}) \bigl(\sigma(y_{t_{j-1}})-
\sigma(x_{0})\bigr) \Delta W_{t_{j}}
\]
and
\[
\mathbf{M}^{(2)}_{T} = \frac{1}{\delta H_{T}} \Biggl( \sum
^{N}_{j=N_{0}} \widetilde{\varkappa}_{j}
\phi_{j} \Delta W_{t_{j}} - \int^{\tau}_{T_{0}}
\widetilde{\phi}(u)\, \mathrm{d}W_{u} \Biggr).
\]
First, we note that the definition
of the sequence $(\widetilde{\varkappa}_{j})_{j\ge1}$ in
\eqref{secSp.4} implies
%
\begin{equation}
\label{secUp.4-1} \sum^{N}_{j=N_{0}} \widetilde{
\varkappa}_{j}\chi_{h,x_{0}}(y_{t_{j-1}}) \le
H_{T} \qquad\mbox{a.s.}
\end{equation}
Therefore, through the condition \eqref{secFc.1}
\begin{eqnarray*}
\mathbf{E}_{\vartheta} \bigl(\mathbf{M}^{(1)}_{T}
\bigr)^{2} &=& \mathbf{E}_{\vartheta} \Biggl( \frac{1}{\delta H^{2}_{T}} \sum
^{N}_{j=N_{0}} \widetilde{\varkappa}^{2}_{j}
\chi_{h,x_{0}}(y_{t_{j-1}}) \bigl( \sigma(y_{t_{j-1}}) -
\sigma(x_{0}) \bigr)^{2} \Biggr)
\\
& \le & \mathbf{E}_{\vartheta} \frac{h^2\sigma^{2}_{\max}}{\delta H_{T}}.
\end{eqnarray*}
Taking into account here that, for $T\ge3$,
%
\begin{equation}
\label{secUp.10-1} H_{T}\ge h(N-N_{0}) (2\sqrt{
\upsilon_{T}}-\upsilon_{T}) \ge h(N-N_{0}) \sqrt{
\upsilon_{T}},
\end{equation}
we obtain
\[
\lim_{T\to\infty} \varphi_{T} \sup_{\vartheta\in\Theta_{\beta}}
\mathbf{E}_{\vartheta} \bigl|\mathbf{M}^{(1)}_{T}\bigr| =0.
\]
Now we study the term $\mathbf{M}^{(2)}_{T}$. To this end, note
that
$t_{\varpi-1}<\tau\le t_{\varpi}$. Therefore,
we can represent this term as
\[
\mathbf{M}^{(2)}_{T}= \frac{1}{\delta H_{T}} \bigl( \varkappa
\phi_{\varpi} \Delta W_{t_{\varpi}} - \phi_{\varpi} (
W_{\tau}- W_{t_{\varpi-1}} ) \bigr).
\]
Moreover, taking into account that the stopping times $\varpi$
and $\tau$ are bounded, one gets
\[
\mathbf{E} (\Delta W_{t_{\varpi}} )^{2}=\delta\quad \mbox{and}\quad
\mathbf{E} (W_{\tau}-W_{t_{\varpi-1}} )^{2} =\mathbf{E} (
\tau-t_{\varpi-1} ) \le\delta.
\]
Therefore, from here and
\eqref{secUp.10-1},
we get
\[
\mathbf{E}_{\vartheta} \bigl(\mathbf{M}^{(2)}_{T}
\bigr)^{2} \le2\mathbf{E}_{\vartheta}\frac{1}{\delta H^2_{T}} \le
\frac{2}{\delta h^2\upsilon_{T}(N-N_{0})^2}
\]
and
\[
\lim_{T\to\infty} \varphi_{T} \sup_{\vartheta\in\Theta_{\beta}}
\mathbf{E}_{\vartheta} \bigl|\mathbf{M}^{(2)}_{T}\bigr| =0.
\]
Due to \eqref{secUp.2} and \eqref{secUp.4}, it is easy to see that
\[
\lim_{T\to\infty} \varphi_{T} \sup_{\vartheta\in\Theta_{\beta}}
\mathbf{E}_{\vartheta} |\Upsilon_{i,T}| =0,\qquad i=1,2.
\]
To put an end to the proof of this theorem, we present the last term on the
right-hand side of \eqref{secUp.7} as
\[
\frac{\sigma(x_{0})}{\sqrt{\delta H_{T}}} \xi_{T}= \frac{\sigma(x_{0})}{\sqrt{\delta h(N-N_{0})}} \biggl(
\frac
{1}{\sqrt{2q_{\vartheta}(x_{0})}}\xi_{T}+K_{T}\xi _{T} \biggr),
\]
where
\[
K_{T}=\frac{1}{\sqrt{2\widetilde{q}_{T}(x_{0})-\upsilon
_{T}}} - \frac{1}{\sqrt{2q_{\vartheta}(x_{0})}},
\]
and we have to show that
%
\begin{equation}
\label{secUp.11}
\lim_{T\to\infty} \sup_{\vartheta\in\Theta_{\beta}}
\mathbf{E}_{\vartheta} \llvert K_{T}\rrvert |\xi_{T}| =0.
\end{equation}
%
It is easy to see that, for any $T>0$, the random variable
$\xi_{T}$ is $(0,1)$-Gaussian conditionally with respect to
$\mathcal{F}_{T_{0}}$.
Therefore,
\[
\mathbf{E}_{\vartheta} \llvert K_{T}\rrvert |\xi_{T}| =
\sqrt{\frac{2}{\uppi}} \mathbf{E}_{\vartheta} \llvert K_{T}
\rrvert.
\]
Taking into account here Lemma~\ref{Le.secSt.3},
we come to the equality \eqref{secUp.11}.
Hence Theorem~\ref{Th.secUp.1}.

\section{Conclusion}

In the paper, we studied the estimation problem of the function $S$ when
its smoothness is known. In the case of unknown smoothness, in order to
construct
an adaptive estimate based on discrete time observations \eqref{secIn.2} in the
model \eqref{secIn.1} we shall use the approach developed in \cite{GaltchoukPergamenshchikov2001} for continuous time observations. The approach
make use of Lepskii's procedure and sequential estimating. Note that
Lepskii's procedure works here just thanks to sequential estimating
since, for the sequential estimate of the function $S$, the stochastic
term in the deviation \eqref{secUp.7}
is a Gaussian random variable. This provides correct estimating the
tail distribution
of a kernel estimate and adapting for the pointwise risk.
Moreover, for adaptive estimating in the case of quadratic risk, we
shall apply the
selection model developed in \cite{GaltchoukPergamenshchikov2011} to
sequential kernel
estimates \eqref{secSp.5}. Note once more, that Gaussianity of the
stochastic term in
\eqref{secUp.7} is a cornerstone result for obtaining a sharp oracle
inequality.
It permits to find Pinsker's constant like to \cite{GaltchoukPergamenshchikov2011}
and then to study the proposed procedure efficiency.

These both programs will be realized in the next paper.

\begin{appendix}
\section*{Appendix}\label{secApp}

\subsection{Geometric ergodicity}
First of all, we recall that in \cite{GaltchoukPergamenshchikov2014}
we have proved the following result.

\begin{theorem}\label{Th.secApp.0}
For any $\varepsilon>0$, there exist constants $R=R(\varepsilon)>0$ and
$\kappa=\kappa(\varepsilon)>0$ such that
\[
\sup_{u\ge0} \mathrm{e}^{\kappa u} \sup_{\|g\|_{*}\le1}
\sup_{x\in\mathbb{R}} \sup_{\vartheta\in\Sigma_{L,M}\times\mathcal{V}} \frac{
|\mathbf{E}_{\vartheta,x}  g(y_{u})
-\mathbf{m}_{\vartheta}(g)|}{
(1+x^2)^{\varepsilon}} \le
R,
\]
where $\mathbf{E}_{\vartheta,x}  (\cdot )=
\mathbf{E}_{\vartheta}  ( \cdot|y_{0}=x )$, $\|
g\|_{*}=\sup_{x}|g(x)|$.
\end{theorem}

\subsection{Concentration inequality}

For any $\mathbb{R}\to\mathbb{R}$ function $f$ belonging to $\mathbf
{L}_{1}(\mathbb{R})$, we set
%
\begin{equation}
\label{secApp.1} \mathbf{D}_{n}(f) = \sum_{k=1}^{n}
\bigl( f(y_{t_{k}}) - \mathbf{m}_{\vartheta}(f) \bigr).
\end{equation}

Now we assume that the frequency $\delta$ in the observations
\eqref{secIn.2} is of the following form
%
\begin{equation}
\label{secApp.1-1}
\delta=\delta_{T}=\frac{1}{(T+1) l_{T}},
\end{equation}
where the function $l_{T}$ is such that,
%
\begin{equation}
\label{secApp.1-2} \lim_{T\to\infty}\frac{l_{T}}{T^{1/2}}=0 \quad\mbox{and}\quad\lim
_{T\to\infty}\frac{l_{T}}{\ln T}=+\infty,
\end{equation}
in particular, the function $l_{T}=(\ln T)^{1+\gamma}$ from
\eqref{secFc.6} is of this kind.
Moreover, let $\varkappa^{*}=\varkappa^{*}_{T}$ be a positive
function satisfying the following properties
%
\begin{equation}
\label{secApp.1-3}
\lim_{T\to\infty} \varkappa^{*}_{T}=0
\quad\mbox{and}\quad \lim_{T\to\infty}\frac{l_{T} (\varkappa^{*}_{T})^{5}}{\ln T}=+\infty.
\end{equation}

\begin{theorem}[(\cite{GaltchoukPergamenshchikov2013})]\label{Th.secApp.1}
Assume that the frequency $\delta$ satisfies \textup{\eqref{secApp.1-1}}--\textup{\eqref{secApp.1-2}}.
Then, for any $a>0$,
%
\begin{equation}
\label{secApp.1-0} \lim_{T\to\infty} T^{a} \sup
_{h\ge T^{-1/2}} \sup_{\vartheta\in\Theta_{\beta}} \mathbf{P}_{\vartheta}
\bigl( \bigl|\mathbf{D}_{N}(\chi_{h,x_{0}})\bigr| \ge \varkappa^{*}_{T}
T \bigr)=0.
\end{equation}
\end{theorem}

\subsection{Proof of Proposition~\texorpdfstring{\protect\ref{Pr.secUp.1}}{4.1}}

First,\vspace*{-3pt} we note that
by the Bunyakovskii--Cauchy--Schwarz inequality
\[
\mathbf{E}_{\vartheta} \bigl( \varrho^{2}_{j}|
\mathcal{F}_{t_{j-1}} \bigr) \le\delta\widetilde{L}^{2} \int
^{t_{j}}_{t_{j-1}} \mathbf{E}_{\vartheta} \bigl(
(y_{u} - y_{t_{j-1}} )^{2} |\mathcal{F}_{t_{j-1}}
\bigr) \,\mathrm{d}u,
\]
where\vspace*{-1.5pt} $\widetilde{L}=\max(L, M)$.
Note now that, for $t_{j-1}\le u\le t_{j}$,
\begin{eqnarray*}
\mathbf{E}_{\vartheta}\bigl( (y_{u} - y_{t_{j-1}}
)^{2} |\mathcal{F}_{t_{j-1}}\bigr) &\le & 2\delta 
 \biggl(\int^{u}_{t_{j-1}} \mathbf{E}_{\vartheta}
\bigl(S^{2}(y_{v})|\mathcal{F}_{t_{j-1}} \bigr)
\,\mathrm{d}v + 
\sigma^{2}_{\max} \biggr)
\\[-2pt]
&\le & 2\delta \biggl(2\int^{u}_{t_{j-1}} \bigl(
M^{2} + L^{2} \mathbf{E}_{\vartheta}\bigl(y^{2}_{v}|
\mathcal{F}_{t_{j-1}} \bigr) \bigr) \,\mathrm{d}v + 
\sigma^{2}_{\max} \biggr).
\end{eqnarray*}
Due to Proposition~\ref{Pr.subsecAMc.7}, we can
estimate the last\vspace*{-1.5pt} conditional expectation as
\[
\sup_{\vartheta\in\Sigma_{L,M}\times[0,\sigma_{\max}]} \sup_{t_{j-1}\le u\le t_{j}}
\mathbf{E}_{\vartheta} \bigl(y^{2}_{u}|
\mathcal{F}_{t_{j-1}} \bigr) \le D_{*}L+y^{2}_{t_{j-1}}
.
\]
Therefore, taking into account that $\chi_{h,x_{0}}(y_{t_{j-1}})
y^{2}_{t_{j-1}}\le x^{2}_{*}$,\vspace*{-1pt} we obtain
%
\begin{equation}
\label{A.2-1} \sup_{t_{j-1}\le u\le t_{j}} \sup_{\vartheta\in\Sigma_{L,M}\times[0,\sigma_{\max}]}
\chi_{h,x_{0}}(y_{t_{j-1}}) \mathbf{E}_{\vartheta} \bigl(
(y_{u} - y_{t_{j-1}} )^{2} | \mathcal{F}_{t_{j-1}}
\bigr) \le L_{1} \delta.
\end{equation}
Therefore,\vspace*{-3pt}
\[
\sup_{j\ge1} \sup_{\vartheta\in\Sigma_{L,M}\times[0,\sigma_{\max}]}
\chi_{h,x_{0}}(y_{t_{j-1}}) \mathbf{E}_{\vartheta} \bigl(
\varrho^{2}_{j}|\mathcal{F}_{t_{j-1}} \bigr) \le
\widetilde{L}^{2} L_{1} \delta^{3}.
\]
Making use of the inequality \eqref{secUp.4-1}
yields the following upper bound, through the Bunyakovskii--Cauchy--Schwarz\vspace*{-2pt} inequality,
\begin{eqnarray*}
\mathbf{E}_{\vartheta} \Upsilon^2_{1,T} &\le&
\mathbf{E}_{\vartheta}\frac{1}{\delta^2 H_{T}} \sum^{N}_{j=N_{0}}
\widetilde{\varkappa}_{j}\chi_{h,x_{0}}(y_{t_{j-1}})
\varrho^2_{j}
\\[-1pt]
&=& \mathbf{E}_{\vartheta}\frac{1}{\delta^2 H_{T}} \sum
^{N}_{j=N_{0}} \widetilde{\varkappa}_{j}
\chi_{h,x_{0}}(y_{t_{j-1}}) \mathbf{E}_{\vartheta} \bigl(
\varrho^2_{j}|\mathcal{F}_{t_{j-1}} \bigr)
\le \widetilde{L}^{2} L_{1} \delta.
\end{eqnarray*}
Hence,\vspace*{-4pt} Proposition~\ref{Pr.secUp.1}.

\subsection{Proof of Proposition~\texorpdfstring{\protect\ref{Pr.secUp.2}}{4.2}}
Note that by the condition\vspace*{-3pt} \eqref{secAe.1}
\begin{eqnarray*}
\mathbf{E}_{\vartheta} \bigl( \bigl(\varrho^{*}_{j}
\bigr)^{2}|\mathcal{F}_{t_{j-1}} \bigr) &=& \int
^{t_{j}}_{t_{j-1}} \mathbf{E}_{\vartheta} \bigl(\bigl(
\sigma(y_{u}) - \sigma(y_{t_{j-1}}) \bigr)^{2} |
\mathcal{F}_{t_{j-1}} \bigr) \,\mathrm{d}u
\\[-2pt]
& \le & \sigma^{2}_{\max} \int^{t_{j}}_{t_{j-1}}
\mathbf{E}_{\vartheta} \bigl( ( y_{u} - y_{t_{j-1}}
)^{2} |\mathcal{F}_{t_{j-1}} \bigr)\, \mathrm{d}u.
\end{eqnarray*}
Therefore, using the inequality \eqref{A.2-1} one has
\[
\sup_{j\ge1} \sup_{\vartheta\in\Sigma_{L,M}\times[0,\sigma_{\max}]}
\chi_{h,x_{0}}(y_{t_{j-1}}) \mathbf{E}_{\vartheta} \bigl( \bigl(
\varrho^{*}_{j}\bigr)^{2}|\mathcal{F}_{t_{j-1}}
\bigr) \le\sigma^{2}_{\max} L_{1}
\delta^{2}.
\]
From here and \eqref{secUp.4-1} and,
taking into account that $0<\widetilde{\varkappa}_{j}\le1$,
we obtain
\begin{eqnarray*}
\mathbf{E}_{\vartheta} \Upsilon^{2}_{2,T} & =&
\mathbf{E}_{\vartheta} \frac{1}{\delta^{2}H^{2}_{T}} \sum^{N}_{j=N_{0}}
\widetilde{\varkappa}^{2}_{j} \chi_{h,x_{0}}(y_{t_{j-1}})
\mathbf{E}_{\vartheta} \bigl( \bigl(\varrho^{*}_{j}
\bigr)^{2} |\mathcal{F}_{t_{j-1}} \bigr)
\\
&\le & \sigma^{2}_{\max} L_{1} \mathbf{E}_{\vartheta}
\frac{1}{H_{T}}.
\end{eqnarray*}
Now the inequality \eqref{secUp.10-1}
yields \eqref{secUp.4}.
Hence, Proposition~\ref{Pr.secUp.2}.

\subsection{Proof of Proposition~\texorpdfstring{\protect\ref{Pr.secUp.3}}{4.3}}

Taking into account the inequalities
\eqref{secUp.4-1} and \eqref{secUp.10-1},
we obtain that, for any $T\ge3$,
\begin{eqnarray*}
\mathbf{E}_{\vartheta} \mathbf{M}^{2}_{T} &=&
\mathbf{E}_{\vartheta} \frac{1}{\delta^{2} H^{2}_{T}} \sum^{N}_{j=N_{0}}
\chi_{h,x_{0}}(y_{t_{j-1}}) \widetilde{\varkappa}^{2}_{j}
\mathbf{E}_{\vartheta} \bigl(\eta^{2}_{j}|
\mathcal{F}_{t_{j-1}} \bigr)
\\
&\le& \mathbf{E}_{\vartheta} \frac{1}{\delta^{2} H^{2}_{T}} \sum
^{N}_{j=N_{0}} \widetilde{\varkappa}_{j}
\chi_{h,x_{0}}(y_{t_{j-1}}) \int^{t_{j}}_{t_{j-1}}
\sigma^{2}(y_{u})\, \mathrm{d}u
\\
&\le& \frac{\sigma^2_{\max}
}{\delta h(N-N_{0})\sqrt{\upsilon_{T}}}.
\end{eqnarray*}
Hence, Proposition~\ref{Pr.secUp.3}.

\subsection{Proof of Proposition~\texorpdfstring{\protect\ref{Pr.secUp.4}}{6.1}}

We start with setting
\[
r_{T}=\frac{(2\widetilde{q}_{T}-\upsilon_{T})h}{\mathbf
{m}_{\vartheta}(\chi_{h,x_{0}})} \quad\mbox{and}\quad N_{1}=N_{0}+r_{T}(N-N_{0}).
\]
Note that
$N_{1}-N_{0}\le (q_{*}\sqrt{\upsilon_{T}}
)^{-1}N:=N^{*}_{1}$,
for sufficiently large $T$.
Moreover, we set
\[
G_{T}= \frac{1}{H_{T}} \sum^{N_{1}}_{j=N_{0}}
f_{h}(y_{t_{j-1}}) \quad\mbox{and}\quad \widehat{G}_{T}=
G_{T}-B_{T}.
\]
Using \eqref{secFc.13}, we can represent the term $G_{T}$
as
\[
G_{T}= 
\frac{N_{1}-N_{0}}{H_{T}} \mathbf{m}_{\vartheta
}(f_{h})
+ 
\frac{1}{H_{T}} \sum^{N_{1}}_{j=N_{0}}
\widetilde{f}_{h}(y_{t_{j-1}}) := G_{1}(T) +
G_{2}(T),
\]
where $\widetilde{f}_{h}(y)=f_{h}(y)-\mathbf{m}_{\vartheta}(f_{h})$.
Taking into account that
$\mathbf{m}_{\vartheta}(\chi_{h,x_{0}})\ge2h q_{*}$, we obtain
\[
\bigl|G_{1}(T)\bigr|= 
\frac{r_{T}(N-N_{0})h}{H_{T}} \bigl|\mathbf{m}^{*}_{\vartheta}(h)\bigr|
\le\frac{1}{2q_{*}}\bigl|\mathbf {m}^{*}_{\vartheta}(h)\bigr|,\qquad
\mathbf{m}^{*}_{\vartheta}(h)= \frac{\mathbf{m}_{\vartheta}(f_{h})}{h}.
\]
Let us represent the last term as
\[
\mathbf{m}^{*}_{\vartheta}(h)=q_{\vartheta}(x_{0})
\Omega _{x_{0},h}(S)+ \widetilde{\mathbf{m}}_{\vartheta}(h),
\]
where $\widetilde{\mathbf{m}}_{\vartheta}(h)
=
\int^{1}_{-1} (
S(x_{0}+hz)
-
S(x_{0})
 )
 (
q_{\vartheta}(x_{0}+hz)
-
q_{\vartheta}(x_{0})
 )
\,\mathrm{d}z$.
Further, by the definition~\eqref{secFc.5}, one has
\[
\Omega_{x_{0},h}(S) = \Omega_{x_{0},h}(S_{0}) +
\Omega_{x_{0},h}(D),
\]
for some function $D$ from $\mathcal{H}^{w}_{x_{0}}(\varepsilon
,\beta,h)$.
Therefore, the properties \eqref{secFc.3}--\eqref{secFc.4} and
\eqref{secFc.6} yield
\[
\lim_{h\to0} \varphi_{T} \sup_{S\in\mathcal{U}_{M}(x_{0},\beta)}
\bigl\llvert \Omega_{x_{0},h}(S) \bigr\rrvert =0.
\]
Obviously, that
\[
\limsup_{h\to0} h^{-2} \sup_{\vartheta\in\Theta_{\beta}} \bigl|
\widetilde{\mathbf{m}}_{\vartheta}(h)\bigr| <\infty.
\]
Hence,
\[
\limsup_{T\to\infty} \varphi_{T} \sup
_{\vartheta\in\Theta_{\beta}} \mathbf{E}_{\vartheta} \bigl|G_{1}(T)\bigr|=0.
\]
Now we note that,
\[
\mathbf{E}_{\vartheta} G^{2}_{2}(T) =
\mathbf{E}_{\vartheta} \frac{1}{H_{T}^2} \Biggl(\sum
^{N_{1}-1}_{j=N_{0}} \Psi_{j} +
\widetilde{f}_{h}^2(y_{t_{N_{1}-1}}) \Biggr),
\]
where
$\Psi_{j}=\widetilde{f}_{h}^2
(y_{t_{j-1}}) + 2\widetilde{f}_{h}(y_{t_{j-1}})
\sum^{N_{1}}_{l=j+1}
\mathbf{E}_{\vartheta} (
\widetilde{f}_{h}
(y_{t_{l-1}})
|\mathcal{F}_{t_{j-1}}
 )$
and $\mathcal{F}_{t}=\sigma\{y_{s}, 0\le s\le t\}$.
Taking into account that
$(y_{t})_{t\ge0}$ is a homogeneous Markov process and that
$|\widetilde{f}_{h}(y)|\le2Mh$, we estimate from above
the last conditional expectation, through
the Theorem~\ref{Th.secApp.0} (for $\varepsilon=1/2$), as
\begin{eqnarray*}
\bigl\llvert \mathbf{E}_{\vartheta} \bigl( \widetilde{f}_{h}
(y_{t_{l-1}}) |\mathcal{F}_{t_{j-1}} \bigr)\bigr\rrvert &=& \bigl
\llvert \mathbf{E}_{\vartheta,y_{t_{j-1}}} \widetilde{f}_{h}
(y_{t_{l-j}}) \bigr\rrvert 
\le2M hR \bigl( 1+y^{2}_{t_{j-1}}
\bigr)^{1/2} \mathrm{e}^{-\kappa t_{l-j}}
\\
&\le & 2M h R \bigl( 1+|y_{t_{j-1}}| \bigr) \mathrm{e}^{-\kappa\delta(l-j)}.
\end{eqnarray*}
Therefore,
\[
|\Psi_{j}|\le 4M^2h^2 \biggl(1+
\frac{2R(1+|y_{t_{j-1}}|)}{\mathrm{e}^{\kappa\delta}-1} \biggr).
\]
From here, bounding $\mathrm{e}^{\kappa\delta}-1$ by $\kappa\delta$, we get
\begin{eqnarray*}
\mathbf{E}_{\vartheta} G^{2}_{2}(T)&\le&
8M^2h^2\mathbf{E}_{\vartheta}\frac{1}{H^2_{T}}\sum
^{N_{1}}_{j=N_{0}} \biggl(1+\frac{R}{\kappa\delta}\bigl(1+|y_{t_{j-1}}|\bigr)
\biggr)
\\
&\le & 8M^2h^2 \biggl(1+\frac{R}{\kappa\delta} \biggr)
\mathbf{E}_{\vartheta}\frac{(N_{1}-N_{0})}{H^2_{T}}
\\
&&{}+ 8M^2h^{2}\frac{R}{\kappa\delta} \mathbf{E}_{\vartheta}
\frac{1}{H^2_{T}} \sum^{N_{1}}_{j=N_{0}} \bigl(
\mathbf{E}_{\vartheta}\bigl(y^2_{t_{j-1}}|\mathcal
{F}_{t_{N_{0}-1}}\bigr) \bigr)^{1/2}.
\end{eqnarray*}
By making use of Proposition~\ref{Pr.subsecAMc.7}, one obtains
\[
\mathbf{E}_{\vartheta} G^{2}_{2}(T)\le
8M^2h^2 \biggl( 1+\frac{R}{\kappa\delta} \biggr)
\mathbf{E}_{\vartheta}\frac{N_{1}-N_{0}}{H^2_{T}} \bigl( 1+\sqrt{D_{*}L}+
|y_{t_{N_{0}-1}}| \bigr).
\]
Now from \eqref{secUp.10-1} it follows that
\[
\frac{N_{1}-N_{0}}{H^2_{T}} = \frac{1}{H_{T} \mathbf{m}_{\vartheta}(\chi_{h,x_{0}})} \le \frac{1}{2h^2\sqrt{\upsilon_{T}}(N-N_{0})q_{*}}.
\]
Thus,
\[
\sup_{\vartheta\in\Theta_{\beta}} \mathbf{E}_{\vartheta} G^{2}_{2}(T)
\le \frac{G^{*}}{\delta\sqrt{\upsilon_{T}}(N-N_{0})},
\]
where $G^{*}=4M^{2}(\kappa+R)
 (1+2\sqrt{D_{*}L}+|y_{0}| )/(\kappa q_{*})$.
From this equality, we obtain immediately
\[
\lim_{T\to\infty} \varphi_{T} \sup_{\vartheta\in\Theta_{\beta}}
\mathbf{E}_{\vartheta} \bigl|G_{2}(T)\bigr| =0.
\]
Let us estimate the term $\widehat{G}_{T}$.
Taking into account the lower bound \eqref{secUp.10-1}, we get
\begin{eqnarray*}
|\widehat{G}_{T}|&\le & \frac{Mh}{H_{T}} \Biggl\llvert \sum
^{N_{1}}_{j=N_{0}} \chi_{h,x_{0}}(y_{t_{j-1}})
- H_{T} \Biggr\rrvert + \frac{2Mh}{H_{T}}
\\
&\le & \frac{M}{\sqrt{\upsilon_{T}}(N-N_{0})} \sum^{N_{1}}_{j=N_{0}}
\bigl\llvert \widetilde{\chi}_{h}(y_{t_{j-1}}) \bigr\rrvert +
\frac{2M}{\sqrt{\upsilon_{T}}(N-N_{0})},
\end{eqnarray*}
where $\widetilde{\chi}_{h}(y)
= \chi_{h,x_{0}}(y)-\mathbf{m}_{\vartheta}(\chi_{h,x_{0}})$.
By making use of the Theorem~\ref{Th.secApp.0} with $\varepsilon=1/2$
one gets
\[
\sum_{j\ge N_{0}} \mathbf{E}_{\vartheta,y_{0}} \bigl|\widetilde{
\chi}_{h}(y_{t_{j-1}})\bigr|\le \frac{R \mathrm{e}^{-\kappa\delta(N_{0}-1)}}{1-\mathrm{e}^{-\kappa\delta}} \bigl(1+
\sqrt{D_{*}L}+|y_{0}| \bigr).
\]
This inequality implies directly
\[
\lim_{T\to\infty} \varphi_{T} \sup_{\vartheta\in\Theta_{\beta}}
\mathbf{E}_{\vartheta} |\widehat{G}_{T}| =0.
\]
Hence, Proposition~\ref{Pr.secUp.4}.

\subsection{Properties of the estimate \texorpdfstring{\protect\eqref{secSp.1-1}}{(3.3)}}

\begin{lemma}\label{Le.secSt.1}
Assume that the parameter $\delta$ is of the
form \eqref{secFc.14}. Then, for any $a>0$,
\[
\lim_{T\to\infty} T^{a} \sup_{\vartheta\in\Theta_{\beta}}
\mathbf{P}_{\vartheta}\bigl(\bigl|\widehat{q}_{T}(x_{0})-q_{\vartheta}(x_{0})\bigr|>
\upsilon_{T}\bigr) =0.
\]
\end{lemma}

\begin{pf}
Denoting $\psi_{\varsigma}(y)=(1/\varsigma)Q((y-x_0)/\varsigma
)$ one has
\begin{eqnarray*}
\widehat{q}_{T}(x_{0})-q_{\vartheta}(x_{0})&=&
\frac{1}{2} \int^{1}_{-1}
\bigl(q_{\vartheta}(x_{0}+\varsigma z)-q_{\vartheta}(x_{0})
\bigr)\, \mathrm{d}z
\\
&&{}+ \frac{1}{2(N_{0}-1)} \mathbf{D}_{N_{0}-1}(\psi_{\varsigma}).
\end{eqnarray*}
Therefore
\begin{eqnarray*}
&& \mathbf{P}_{\vartheta}\bigl(\bigl|\widehat{q}_{T}(x_{0})-q_{\vartheta}(x_{0})\bigr|>
\upsilon_{T}\bigr)
\\
&& \quad \le \mathbf{P}_{\vartheta}\biggl(\biggl|\int^{1}_{-1}
\bigl(q_{\vartheta
}(x_{0}+\varsigma z)-q_{\vartheta}(x_{0})
\bigr) \,\mathrm{d}z\biggr|> \upsilon_{T}\biggr) +\mathbf{P}_{\vartheta}
\biggl(\frac{1}{N_{0}-1} \mathbf{D}_{N_{0}-1}(\psi_{\varsigma})>
\upsilon_{T}\biggr).
\end{eqnarray*}
The first term on the right-hand side equals to zero for sufficiently
large $T$ since
\[
\biggl|\int^{1}_{-1} \bigl(q_{\vartheta}(x_{0}+
\varsigma z)-q_{\vartheta}(x_{0})\bigr) \,\mathrm{d}z\biggr| \le
\varsigma^2\ddot{q}^*<\upsilon_{T},
\]
for sufficiently large $T$, where $\ddot{q}^*=\sup_x\sup_{\vartheta
}|\ddot{q}_{\vartheta}(x)|<\infty$.
Applying Theorem~\ref{Th.secApp.1}, to the second term on the
right-hand side of the
same inequality yields the Lemma~\ref{Le.secSt.1}.
\end{pf}

\begin{lemma}\label{Le.secSt.2}
Assume that the parameter $\delta$ is of the
form \eqref{secFc.14}. Then, for any $a>0$,
\[
\lim_{T\to\infty} T^{a} \sup_{\vartheta\in\Theta_{\beta}}
\mathbf{P}_{\vartheta}\bigl(\bigl|\widetilde{q}_{T}(x_{0})-q_{\vartheta}(x_{0})\bigr|>
\upsilon_{T}\bigr) =0.
\]
\end{lemma}

\begin{pf}
Note, that for sufficiently large $T$
(for that $\ln T \ge\max(q^{*2}, 1/q^{2}_{*})$),
\[
\bigl|\widetilde{q}_{T}(x_{0})-q_{\vartheta}(x_{0})\bigr|
\le \bigl|\widehat{q}_{T}(x_{0})-q_{\vartheta}(x_{0})\bigr|.
\]
The lemma follows immediately from Lemma~\ref{Le.secSt.1}.
\end{pf}

%
\begin{lemma}\label{Le.secSt.3}
Assume that the parameter $\delta$ is of the
form \eqref{secFc.14}. Then,
%
\begin{equation}
\label{A.3} \limsup_{T\to\infty} \frac{1}{\upsilon^{1/2}_{T}} \sup
_{\vartheta\in\Theta_{\beta}} \mathbf{E}_{\vartheta} \biggl\llvert
\frac{1}{\widetilde{q}_{T}(x_{0})-\upsilon_{T}} - \frac{1}{q_{\vartheta}(x_{0})} \biggr\rrvert \le\frac{4}{q^*} <
\infty.
\end{equation}
\end{lemma}

\begin{pf}
Indeed, for sufficiently large $T$
for which
\[
\upsilon_{T} \le \min \bigl( 1/\bigl(q^{*}
\bigr)^{2}, 1/4 \bigr),
\]
we obtain
\begin{eqnarray*}
\mathbf{E}_{\vartheta} \biggl\llvert \frac{1}{\widetilde{q}_{T}(x_{0})-\upsilon_{T}} - \frac{1}{q_{\vartheta}(x_{0})}
\biggr\rrvert & \le& \frac{2\upsilon^{1/2}_{T}}{q_{*}} + \frac{2}{q_{*}\upsilon^{1/2}_{T}} \mathbf{E}_{\vartheta}
\bigl\llvert \widetilde{q}_{T}(x_{0}) -
q_{\vartheta}(x_{0}) \bigr\rrvert
\\
&\le& \frac{4\upsilon^{1/2}_{T}}{q_{*}} + \frac{2}{q_{*}\upsilon_{T}} \mathbf{P}_{\vartheta} \bigl( \bigl|
\widetilde{q}_{T}(x_{0}) - q_{\vartheta}(x_{0})
\bigr|>\upsilon_{T} \bigr).
\end{eqnarray*}
Now Lemma~\ref{Le.secSt.2} implies the equality \eqref{A.3}.
Hence, Lemma~\ref{Le.secSt.3}.
\end{pf}

\subsection{Moment inequality for the process \texorpdfstring{\protect\eqref{secIn.1}}{(1.1)}}

We state the moment bound from \cite{GaltchoukPergamenshchikov2014}.

\begin{proposition}\label{Pr.subsecAMc.7}
Let $(y_{t})_{t\ge0}$ be a solution of the equation \eqref
{secIn.1}. Then,
for any $z\in\mathbb{R}$ and $m\ge1$,
%
\begin{equation}
\label{subsecAMc.15} \sup_{u\ge0} \sup_{\vartheta\in\Sigma_{L,M}\times[0,\sigma_{\max}]}
\mathbf{E}_{\vartheta,z} (y_{u})^{2m} \le(2m-1)!!
\bigl(D_{*}L+ z^{2}\bigr)^{m},
\end{equation}
where $\mathbf{E}_{\vartheta,z}(\cdot)=\mathbf{E}_{\vartheta}(\cdot|y_{0}=z)$ and
$D_{*}=(
M+Lx_{*}+2x_{*}
)^{2}
(L+M)+\sigma^{2}_{\max}$.
\end{proposition}

\begin{pf}To obtain this inequality, we make use of the method proposed
in \cite{KabanovPergamenshchikov2003}, page 20,  for linear stochastic
equations. First of
all, note that thanks to Theorem~4.7 from \cite{LiptserShiryaev1978},
for any
$T>0$, there exists some $\varepsilon>0$ such that, for each
$\vartheta\in\Theta_{\beta}$ and $z\in\mathbb{R}$,
%
\begin{equation}
\label{subsecAMc.15-1} \sup_{0\le t\le T} \mathbf{E}_{\vartheta,z}
\mathrm{e}^{\varepsilon y^{2}_{t}} <\infty.
\end{equation}
Let us denote
$D_{\vartheta}(y)=2yS(y)+\sigma^{2}(y)+\check{L} y^{2}$
and $\check{L}=L^{-1}$.
Taking into account that $0<\check{L}<1$ and $x_{*}\ge1$,
we obtain that, for $|y|\le x_{*}$,
\[
\bigl|D_{\vartheta}(y)\bigr|\le x^{2}_{*}(2 M+1) +
\sigma^{2}_{\max}.
\]
Let now $|y|\ge x_{*}$. Denoting by $y_{*}$ the projection of
$y$ onto
the interval $[-x_{*}, x_{*}]$ we obtain that
\begin{eqnarray*}
2yS(y)&=& 2yS(y_{*})+2 y_{*} \bigl(S(y)-S(y_{*})
\bigr) +2(y-y_{*}) \bigl(S(y)-S(y_{*}) \bigr)
\\
&\le &  2|y| M + 2L x_{*} |y-y_{*}| -2\check{L}
|y-y_{*}|^{2}
\\
&\le & 2(M+Lx_{*}+2x_{*})|y| -2\check{L} y^{2}.
\end{eqnarray*}
Therefore,
\[
\sup_{\vartheta\in\Sigma_{L,M}\times[0,\sigma_{\max}]} \sup_{y\in\mathbb{R}} D_{\vartheta}(y)
\le D_{*}.
\]
By
the It\^{o} formula, we obtain
\begin{eqnarray*}
\mathrm{d}y_{u}^{2m}&=& -m\check{L} y_{u}^{2m}
\,\mathrm{d}t+ m y_{u}^{2(m-1)} \bigl(D_{\vartheta}(y_{u})+2(m-1)
\sigma ^{2}(y_{u}) \bigr)\,\mathrm{d}t
\\
&&{}+2m y_{u}^{2m-1} \sigma(y_{u})
\,\mathrm{d}W_{t}.
\end{eqnarray*}
Moreover, the property \eqref{subsecAMc.15-1} yields that, for any
$m\ge1$,
\[
\mathbf{E}_{\vartheta} \int^{t}_{0}
\mathrm{e}^{-m\check{L}(t-s)} y_{s}^{2m-1} \sigma (y_{s})
\,\mathrm{d}W_{s}=0.
\]
Therefore, $\mathbf{E}_{\vartheta} y^{2m}_{t}
\le z^{2m}+m(2m-1)D_{*}
 \int^{t}_{0} \mathrm{e}^{-m\check{L} (t-s)}
\mathbf{E}_{\vartheta} y^{2(m-1)}_{s}
\,\mathrm{d}s$.
Now the induction implies directly the bound \eqref{subsecAMc.15}.
Hence, Proposition~\ref{Pr.subsecAMc.7}.
\end{pf}

\begin{proposition}\label{Pr.subsecAMc.8}
Let $(y_{t})_{t\ge0}$ be a solution of the equation \eqref{secIn.1}.
Then, for any $z\in\mathbb{R}$ and $m\ge1$, and for any
stopping time $\tau$ taking values in $[0,T]$, one has
%
\begin{equation}
\label{subsecAMc.16} \sup_{\vartheta\in\Sigma_{L,M}\times[0,\sigma_{\max}]} \mathbf{E}_{\vartheta,z}
(y_{\tau})^{2m} \le B^{*}(m,z) T
\end{equation}
and
%
\begin{equation}
\label{subsecAMc.17} \sup_{\vartheta\in\Sigma_{L,M}\times[0,\sigma_{\max}]} \mathbf{E}_{\vartheta,z} \sup
_{0\le u\le T}(y_{u})^{2m} \le
B^{*}_{1}(m,z) T,
\end{equation}
where
$B^{*}(m,z)=(2m-1)!!
(D_{*}L+z^{2})^{m}
(D_{*}+2(m-1)\sigma^{2}_{\max})$
and $B^{*}_{1}(m,z)=1+mB^{*}(m+1,z)$.
\end{proposition}

The proof of this proposition follows immediately from Proposition~1.1.5
in \cite{KabanovPergamenshchikov2003}.

\subsection{Proof of Proposition~\texorpdfstring{\protect\ref{Pr.secNa.1}}{4.5}}

It is clear, that to show \eqref{secSp.5} it suffices to check that,
for any $a>0$,
%
\begin{equation}
\label{A.2} \lim_{T\to\infty} T^{a} \sup
_{\vartheta\in\Sigma_{L,M}\times[\sigma_{\min},
\sigma_{\max}]} \mathbf{P}_{\vartheta}\bigl(\Gamma^{c}_{T}
\bigr)=0.
\end{equation}
Indeed, by the definition of $\varpi$
\begin{eqnarray*}
\mathbf{P}_{\vartheta}\bigl(\Gamma^{c}_{T}\bigr)&=&
\mathbf{P}_{\vartheta} \Biggl(\sum^{N}_{j=N_{0}}
\chi_{h,x_{0}}(y_{t_{j-1}}) <H_{T} \Biggr)
\\
&=&\mathbf{P}_{\vartheta} \bigl(\mathbf{D}_{N_{0},N-1}(\chi
_{h,x_{0}}) < \bigl(2\widetilde{q}_{T}-\upsilon_{T}-
\mathbf{m}^{*}_{\vartheta}(\chi_{h,x_{0}})\bigr)
(N-N_{0})h \bigr),
\end{eqnarray*}
where $\mathbf{D}_{k,n}(f)=\mathbf{D}_{n}(f) -\mathbf{D}_{k}(f)$ and
\[
\mathbf{m}^{*}_{\vartheta}(\chi_{h,x_{0}})=\frac
{\mathbf{m}_{\vartheta}(\chi_{h,x_{0}})}{h}
=\int^{1}_{-1} q_{\vartheta}(x_{0}+hz)
\,\mathrm{d}z.
\]
Taking into account the definition of $\upsilon_{T}$ in \eqref{secSp.1-1}
we obtain that, for sufficiently large $T$,
\[
\sup_{\vartheta\in\Theta_{\beta}} \int^{1}_{-1}
\bigl|q_{\vartheta}(x_{0}+hz)-q_{\vartheta}(x_{0})\bigr|
\,\mathrm{d}z\le\upsilon_{T}/4.
\]
Therefore, for such $T$,
\begin{eqnarray*}
\mathbf{P}_{\vartheta}(\varpi>N) &\le& \mathbf{P}_{\vartheta} \bigl(\bigl|
\widetilde{q}_{T}(x_{0})- q_{\vartheta}(x_{0})\bigr|>
\upsilon_{T}/8 \bigr)
\\
&&{}+ \mathbf{P}_{\vartheta} \bigl( \bigl|\mathbf{D}_{N_{0},N-1}(
\chi_{h,x_{0}})\bigr| > N h\upsilon_{T}/2 \bigr).
\end{eqnarray*}
Now we estimate the last term as
\begin{eqnarray*}
\mathbf{P}_{\vartheta} \bigl( \bigl|\mathbf{D}_{N_{0},N-1}(
\chi_{h,x_{0}})\bigr| > N h\upsilon_{T}/2 \bigr)&\le&
\mathbf{P}_{\vartheta} \bigl( \bigl|\mathbf{D}_{N-1}(\chi_{h,x_{0}})\bigr|
> N h\upsilon_{T}/4 \bigr)
\\
&&{}+ \mathbf{P}_{\vartheta} \bigl( \bigl|\mathbf{D}_{N_{0}}(
\chi_{h,x_{0}})\bigr| > N h\upsilon_{T}/4 \bigr).
\end{eqnarray*}
By applying Lemma~\ref{Le.secSt.1} and the inequality
\eqref{secApp.1-0}, we obtain (\ref{A.2}).
\end{appendix}

\section*{Acknowledgements}
The authors are grateful to the Associated Editor and to the anonymous
referee for a careful reading of manuscript
and helpful comments which led to an improved presentation of the
paper.

The second author is  supported in part by
the Russian Science Foundation (the research project No. 14-49-00079)
and
the Mathematics and Mechanics Department of the National Research Tomsk
State University,
Tomsk, Russia.


\printhistory
\end{document}